\documentclass[a4paper]{article}

\usepackage{amsmath,amssymb}
\usepackage{chngcntr,theorem,mathdots}

\counterwithin{equation}{section}

\newtheorem{theorem}[equation]{Theorem}

\newtheorem{lemma}[equation]{Lemma}

\newtheorem{cor}[equation]{Corollary}

{\theorembodyfont{\rmfamily}
 }
{\theorembodyfont{\rmfamily}
 }
{\theorembodyfont{\rmfamily}
 }
{\theorembodyfont{\rmfamily}
 \newtheorem{defi}[equation]{Definition}}

\newenvironment{romanlist}
  {%
   \setlength{\topsep}{0pt}%
   \vspace{-\parskip}%
   \begin{enumerate}%
     \setlength{\parsep}{0pt}
     \setlength{\parskip}{0pt}%
  }%
  {\end{enumerate}%
   \vspace{-\parskip}}

\newcommand{\proofend}{\hfill $ \framebox{} $ \par}

\newcommand{\Z}{\mathbb{Z}} \newcommand{\N}{\mathbb{N}}

\newcommand{\dotcup}{\ensuremath{\,\,\mathaccent\cdot\cup\,\,}}
\newcommand{\boxrel}{\ensuremath{\,\,\Box\,\,}}
\newcommand{\strgpr}{\ensuremath{\,\boxtimes\,}}

\title{Factors of disconnected graphs and polynomials with nonnegative
integer coefficients}
\author{Christiaan van de Woestijne\\
Lehrstuhl f\"ur Mathematik und Statistik\\
Montanuniversit\"at Leoben\\
Franz-Josef-Stra\ss e 18, 8700 Leoben, Austria\\
{\tt c.vandewoestijne@unileoben.ac.at}}
\date{20 December 2010}

\begin{document}

\maketitle

\begin{abstract}
  We investigate the uniqueness of factorisation of possibly disconnected
  finite graphs with respect to the Cartesian, the strong and the direct
  product. It is proved that if a graph has $n$ connected components, 
  where $n$ is prime, or $n=1,4,8,9$, and satisfies some additional
  conditions, it factors uniquely under the given products. If, on the
  contrary, $n=6$ or $10$, all cases of nonunique factorisation are described
  precisely.
\end{abstract}

\section{Introduction}

The purpose of this note is to give an algebraic/combinatoric perspective on a
result in \cite{HammackImrichKlavzar:2011} about uniqueness of factorisation of
certain finite or infinite graphs with respect to well-known graph products. As
already noted by the authors of \cite{HammackImrichKlavzar:2011}, the structure
underlying the problem is the factorisation of multivariate polynomials with
nonnegative integer coefficients. We expand in this theme, and provide precise
classifications of all possible non-unique factorisations for a number of
small cases, which include polynomials of unbounded degree, however. These
results immediately translate back to graphs, and as our main result, we have
the following (definitions are given below).

\begin{theorem}
  Let $G$ be a graph having $n$ connected components, where either $n=1$,
  $4$, $8$, or $9$, or $n$ is a prime number. Then if $G$ is simple, it has
  unique Cartesian and strong factorisation; if $G$ has no bipartite
  components, it has unique direct factorisation. If, on the contrary,
  $n=6$, $10$, $12$, or $16$, non-unique factorisations occur for all three
  products.
\end{theorem}

\paragraph*{Proof.}
  By the Structure Theorem \ref{ThmStruc}, together with Lemma \ref{LemNovars},
  we reduce the problem to the unique or non-unique factorisation of univariate
  polynomials $P$ with nonnegative integer coefficients such that $P(1)=n$.
  This problem is trivial for $n=1$, easy for $n$ prime (see Lemma
  \ref{LemPrime}), and solvable with the aid of a computer for the cases
  $n\in \{4,6,8,9,10\}$ (see Theorems \ref{Thm4} up to \ref{Thm10}); for
  $n=12$ and $16$, we have the examples \eqref{EqNusken} and \eqref{Eq16}.
\proofend

The results for $n$ prime, or $n=4$, are given as Exercises 6.12 and 6.13 in
\cite{HammackImrichKlavzar:2011}. The authors provide directly
graph-theoretical proofs, so that we obtain a second proof for unique
factorisation of polynomials with nonnegative integer coefficients and $4$ or a
prime number of terms.

The problem of unique factorisation of (disconnected or even connected) graphs
\emph{with loops} with respect to the Cartesian and the strong product is still
open --- cf.\ Lemma \ref{LemLoop} below.

\section{Factorisations of graphs}

We will consider the class of finite graphs (or symmetric binary relations)
$\Gamma_0$ and its subclass consisting of graphs without loops (or simple
graphs) $\Gamma$. (To be precise, we consider the category of \emph{isomorphism
classes} of such graphs, with graph homomorphisms as arrows. This will be
tacitly assumed in the sequel.) A graph $G$ is given as a pair $(V,E)$, where
$V$ is the set of vertices and $E\subseteq V\times V$ is the set of edges. Note
that $(x,y)\in E \Leftrightarrow (y,x)\in E$, as our graphs are undirected.

\begin{defi}
  Let $G$ and $H$ be graphs in $\Gamma_0$.
  \begin{romanlist}
    \item
      The \emph{disjoint union} $G\dotcup H$ is the graph with vertex set $V(G)
      \dotcup V(H)$ and edge set $E(G) \dotcup E(H)$ (both again disjoint
      union).
    \item
      The \emph{Cartesian product} $G\boxrel H$ is the graph with vertex set
      $V(G)\times V(H)$, such that $(u,v)$ and $(u',v')$ are adjacent if and
      only if $u=u'$ and $(v,v')\in E(H)$ or $(u,u')\in E(G)$ and $v=v'$.
    \item
      The \emph{strong product} $G \strgpr H$ is the graph with vertex set
      $V(G)\times V(H)$, such that distinct vertices $(u,v)$ and $(u',v')$ are
      adjacent if and only if $u=u'$ or $(u,u')\in E(G)$, and $v=v'$ or
      $(v,v')\in E(H)$, and $(u,v)$ has a loop if and only if at least one of
      $u$ and $v$ has a loop.
    \item
      The \emph{direct product} $G\times H$ is the graph with vertex set
      $V(G)\times V(H)$,\\
      such that $(u,v)$ and $(u',v')$ are adjacent if and
      only if $(u,u')\in E(G)$ and $(v,v')\in E(H)$.
  \end{romanlist}
\end{defi}

All three products are commutative and associative and distribute over the
disjoint union, and are the only products with these properties such that the
projections onto the factors are graph homomorphisms (direct product) or weak
homomorphisms (the other two). Details can be found in
\cite{HammackImrichKlavzar:2011}, who only define the Cartesian and strong
products for simple graphs, but the extension to arbitrary graphs, using the
above definitions, is immediate. The graph $K_1$ consisting of one vertex and
no edges is a neutral element with respect to both the Cartesian and the strong
product, whereas the graph $K_1^*$, having one vertex and a loop on it, is a
neutral element for the direct product. Of course, the empty graph is a neutral
element for the disjoint union. Finally, note that the Cartesian and strong
products of two graphs are simple if and only if both factors are simple. It
follows that both $\Gamma$ and $\Gamma_0$ are commutative semirings with
operations given by the disjoint union and either the Cartesian or the strong
product, whereas $\Gamma_0$ is a commutative semiring with the disjoint union
and the direct product.

An essential question regarding graph products is whether the possible
factorisations obtained with respect to such a product are unique. Here, we
have the following results.

\begin{defi}
  Let $\odot$ be a graph product. A graph $G$, not the neutral element for
  $\odot$, is \emph{irreducible} with respect to $\odot$ if in any
  factorisation $G=H\odot L$, either $H$ or $L$ is the neutral element for
  $\odot$.
\end{defi}

\begin{theorem} \cite[Theorems 6.6, 7.14, 8.17]{HammackImrichKlavzar:2011}\\%
  \begin{romanlist}
    \topsep 0pt%
    \vspace{-\baselineskip}%
    \item
      (Sabidussi, Vizing) Let $G$ be a connected simple graph. Then any two
      factorisations of $G$ into irreducibles with respect to the Cartesian
      product are unique up to permutation and isomorphism of the factors.
    \item
      (D\"orfler and Imrich, McKenzie) The same result holds for the strong product.
    \item
      (McKenzie) Let $G$ be a connected nonbipartite graph. Then any two
      factorisations of $G$ into irreducibles with respect to the direct
      product are unique up to permutation and isomorphism of the factors.
  \end{romanlist}
\end{theorem}

The next result shows that Cartesian and strong factorisation for graphs with
loops is in general not unique. The proof is easy.

\begin{lemma} \label{LemLoop}
  Let $G$ be a graph in which every node has a loop. Then the same holds for
  $G\boxrel H$ and $G\strgpr H$, where $H$ is any graph. 
\end{lemma}

Bipartite graphs, with or without loops, generally have several nonequivalent
direct factorisations. A classification of these phenomena is a topic of
current research by Hammack and his co-authors
\cite{Abay-AsmeromHammackLarsonTaylor:2009/2010,
Hammack:2006:Isom,Hammack:2009:Proo}.

In this note, we will consider factorisations of \emph{disconnected} graphs. In
this case, there are certain well-known examples of non-unique factorisations.
However, we will show that these examples are all reducible to the same
phenomenon in the semiring of polynomials in several variables with nonnegative
integral coefficients. 

To do this, we slightly change our perspective on the class of graphs. Instead
of thinking of the connected components of a graph as being separate graphs
standing next to each other, and possibly repeated more than once, we will
think of each connected component as occurring only once, but with some
\emph{multiplicity} assigned to it. Then, an easy step is to allow also
negative multiplicities.

\begin{defi}
  We define the class of \emph{graphs with integer-weighted components}
  $\tilde{\Gamma}_0$ to be the set of formal sums $\sum_G a_G G$, where $G$
  runs over the \emph{connected} graphs in $\Gamma_0$, and where the $a_G$ are
  integers, only finitely many of which are nonzero. We define
  \begin{align*}
    (\sum_G a_G G) + (\sum_G b_G G) &= \sum_G (a_G + b_G) G; \\
    (\sum_G a_G G) \odot (\sum_G b_G G) &= 
      \sum_G \left( \sum_{H\odot L = G} a_H b_L \right) G.
  \end{align*}
  Doing exactly the same for the subclass $\Gamma$ of simple graphs yields
  the set $\tilde{\Gamma}$.
\end{defi}

In the definition of multiplication, we obtain a double sum, because it is
possible that several combinations of components $H$ on the left and $L$ on the
right have the same product $H\odot L=G$. We combine all these and add the
multiplicities together to obtain the multiplicity of (the isomorphism class
of) $G$ in the product.

We note that we will continue to use the term \emph{graph} in the usual way;
the term \emph{graph with integer-weighted components} will denote a formal sum
of connected graphs, with possibly negative multiplicities.

\begin{theorem}
  For any choice of graph product $\odot=\boxrel$, $\strgpr$, $\times$, the set
  $\tilde{\Gamma}_0$ is a commutative ring with the addition and multiplication
  as defined above, and the semiring $\Gamma_0$ embeds into it. The additive
  group of $\tilde{\Gamma}_0$ is the free Abelian group generated by the 
  connected graphs in $\Gamma_0$. 
  
  Analogously, $\tilde{\Gamma}$ is a ring, and the semiring $\Gamma$ embeds
  into it. The connected simple graphs freely generate its additive group.
\end{theorem}

\paragraph*{Proof.}
  It is clear from the definition that the set $\tilde{\Gamma}_0$ is just the
  free Abelian group generated by the connected graphs.

  If in a formal sum all weights $a_G$ are \emph{nonnegative}, then the formal
  sum is just a graph as we know it. Conversely, every graph can be written as
  such a sum in exactly one way, because the decomposition of graphs as a
  disjoint union of connected graphs is unique. Thus, we can embed the set of
  graphs into the set $\tilde{\Gamma}_0$.

  Furthermore, one sees that addition and multiplication of nonnegatively
  weighted graphs agrees exactly with the operations of $\dotcup$ and $\odot$
  on graphs, so that the embedding is a homomorphism of semirings.

  Now what remains is to consider the operations on the whole set
  $\tilde{\Gamma}_0$. The operations of addition and multiplication are clearly
  commutative. Their associativity and the distributivity of $\odot$ over
  addition are easily derived from the corresponding properties of
  nonnegatively weighted graphs, using the \emph{difference representation} of
  formal sums, as follows. Every graph with integer-weighted components $X$ can
  be written as a difference $X_1 - X_2$, where both $X_1$ and $X_2$ have only
  nonnegative multiplicities.  If we also say that every connected graph has
  nonzero multiplicity in at most one of $X_1$ and $X_2$, this decomposition is
  uniquely determined, and we see that elements of $\tilde{\Gamma}_0$ can be
  uniquely written as the \emph{difference} $G-H$ of two graphs $G$ and $H$
  that have no connected component in common (i.e., no component of $G$ is
  isomorphic to any component of $H$ and vice versa). The details are left to
  the reader.
\proofend

~

The construction used here is known as forming the \emph{Grothendieck group} of
a commutative semigroup with cancellation. When the semigroup is a semiring (as
in our case), the proof shows that the construction easily extends to define a
compatible ring structure on the Grothendieck group.

Now we come to the main structure theorem (see also \cite{FernandezLeightonLopez-Presa:2007, ImrichKlavzarRall:2007}).

\begin{theorem} \label{ThmStruc}
  The following three rings are each isomorphic to the ring
  $R=\Z[x_1,x_2,\ldots, x_n,\ldots]$ of polynomials with integer coefficients
  in countably many commuting variables:
  \begin{romanlist}
    \item
      the ring $\tilde{\Gamma}$ of simple graphs with the Cartesian product;
    \item
      the ring $\tilde{\Gamma}$ of simple graphs with the strong product;
    \item
      the subring $\tilde{\Gamma}_{0,\rm odd}$ of $\tilde{\Gamma}_0$ of graphs
      having only nonbipartite connected components, with the direct product.
  \end{romanlist}
\end{theorem}

Note that, in fact, we have unique direct factorisation for all nonbipartite
graphs \cite[Theorem 8.17]{HammackImrichKlavzar:2011}. However, the set of all
these is not a free Abelian group under disjoint union: if $G$ is connected
bipartite and $H$, $L$ are connected, nonbipartite, and nonisomorphic, then we
have
$$
  (G\dotcup H) \dotcup L = (G\dotcup L) \dotcup H
$$
where none of the factors can be written as a sum of nonempty graphs within the
class of nonbipartite graphs. One might call this kind of phenomenon a
``non-unique partition''.

\paragraph*{Proof.}
  Let $\odot$ be either the Cartesian or the strong product. As unique
  factorisation is only known for connected simple graphs, we restrict
  ourselves to the ring $\tilde{\Gamma}$.

  Now let $G_1,G_2,\ldots$ be some (arbitrary) enumeration of the set of
  \emph{irreducible} connected simple graphs.

  We now define a ring homomorphism $\phi:R\rightarrow \tilde{\Gamma}$ by
  sending the variable $x_i$ to $G_i$ (``evaluating the polynomial at the
  values $G_i$, for $i\ge 1$''); it is obvious that this determines the
  homomorphism, by linearity. Of course, the unit element $1$ is sent to the
  neutral element $K_1$, and $0$ is sent to the empty graph. As every connected
  simple graph has a unique factorisation into irreducible graphs, the images
  of the monomials $x_{i_1}^{e_1}\cdot \ldots \cdot x_{i_r}^{e_r}$ contained in
  $R$ are all distinct, and every connected simple graphs is the image of one
  such monomial. Thus, because the generators of the additive group of $R$
  (which are the monomials) are mapped one-to-one onto the generators of the
  additive group of $\tilde{\Gamma}$, and both groups are free Abelian, we find
  that $\phi$ is an isomorphism.

  Next, consider the direct product. It is an easy exercise to show that the
  direct product of two connected graphs in $\Gamma_0$ is bipartite if and only
  if at least one of the factors is bipartite. This means that the subgroup 
  $\tilde{\Gamma}_{0,\rm odd}$ of $\tilde{\Gamma}_0$ generated by the
  nonbipartite connected graphs is closed under taking direct products. As it
  is obviously a subgroup with respect to the disjoint union and contains the
  neutral element $K_1^*$, it is a subring.

  Now exactly the same argument as before shows that this ring is isomorphic to
  $R$.
\proofend

~

Of the many consequences of the above isomorphisms, we only give the following.
Recall that a set $S$ has the \emph{cancellation property} with respect to some 
binary operation $\odot$ on $S$, if for all $G,H,L\in S$, the isomorphism
$G\odot H \cong G\odot L$ implies the isomorphism $H\cong L$. It means that
whenever we find that a certain element of $S$ has a factor $G$ (is
``divisible'' by $G$), the ``quotient'' by $G$ is uniquely determined within
$S$. Note that the cancellation property is inherited by any subset of $S$ that
is closed under the operation $\odot$, whereas unique factorisation does not
necessarily descend in this way.

The following result is due to Lov\'asz \cite[Theorem 9]{Lovasz:1971:Canc} for
the case of the direct product; see also \cite[Theorem
9.10]{HammackImrichKlavzar:2011}. In fact, what Lov\'asz proved was that
quotients by nonbipartite simple graphs and by reflexive nonbipartite graphs
are well-defined within the class of (not necessarily symmetric) binary
relations. However, his argument extends almost verbatim to the case of
nonbipartite graphs that may have some loops, and the result follows by the
inheritance property just mentioned. Of course, our proof only shows that the
cancellation law follows directly from the unique factorisation of connected
graphs. As said above, we restrict to graphs without bipartite components, to
avoid problems with non-unique partition.

For the Cartesian and strong cases, the argument appears as \cite[Theorems 6.21
and 9.5]{HammackImrichKlavzar:2011}.

\begin{cor}
  The set of nonempty simple graphs has the cancellation property with respect
  to the Cartesian and the strong product. The set of graphs without bipartite
  components has the cancellation property with respect to the direct product.
\end{cor}

\paragraph*{Proof.}
  The given sets are embedded (with the operation being respected) in the ring
  $R$, which is a unique factorisation domain and hence certainly has the
  cancellation property with respect to multiplication of polynomials.

  As any two distinct factorisations $G\odot H\cong G\odot L$ would lead to
  distinct factorisations of the same polynomial, we conclude that cancellation
  holds.
\proofend

\section{Polynomials with nonnegative integer coefficients}

We now proceed to the actual purpose of this note. By the above isomorphisms,
it is easy to exhibit examples of disconnected graphs that have multiple
non-equivalent factorisations. For example, consider the polynomial identity
\begin{equation} \label{EqNonunique}
  (1+x+x^2)(1+x^3) = (1+x^2+x^4)(1+x).
\end{equation}
We have factorisations $1+x^3=(1+x)(1-x+x^2)$ and
$1+x^2+x^4=(1-x+x^2)(1+x+x^2)$, but as these contain factors with negative
coefficients, they do not correspond to honest graph factorisation. Restricting
ourselves to the semiring $N$ of polynomials with nonnegative integer
coefficients, it follows that all the factors occurring in \eqref{EqNonunique}
are irreducible, and hence so are the corresponding graphs, if we let $x$
correspond to some irreducible connected simple graph of our choice (as $K_1$
is a ``unit'', we do not consider it to be irreducible). Mapping $x$ to $K_2$,
the graph with two nodes and one edge, we obtain
$$
  (K_1 + K_2 + K_2^{\odot,2}) \odot (K_1 + K_2^{\odot,3}) = 
  (K_1 + K_2^{\odot,2} + K_2^{\odot,4}) \odot (K_1 +K_2),
$$
where $\odot$ is either the Cartesian or the strong product, and
$G^{\odot,n}$ means the $n$-fold $\odot$-product of $G$ with itself (also
called the $n$th $\odot$-power of $G$).

For the direct product, we can use the same example. However, here we must use
$K_1^*$ (the vertex with a loop) as the neutral element, and the choice for
$\phi(x)$ must be some irreducible connected \emph{nonbipartite} graph, such as
$K_2^*$ (two vertices with loops and an edge between them).

Another example, due to N\"usken \cite{FernandezLeightonLopez-Presa:2007} and
containing nontrivial integer coefficients, is
\begin{equation} \label{EqNusken}
  (2+x+x^3)(2+x)=(1+x)(4+x^2+x^3).
\end{equation}
Here the factor $2-x+x^2$ divides both cubic factors, and it follows that
the latter are irreducible in the semiring $N$.

However, by the structure theorem \ref{ThmStruc}, it follows that 
\emph{every instance of non-unique factorisation} existing among graphs arises
from a non-unique factorisation in $N$. We can use this fact to limit the
possible phenomena inside our sets of graphs. Every graph corresponds to
a polynomial, with each connected component corresponding to a monomial, whose
coefficient is equal to the number of components that share the same
isomorphism class.

Hence, what we will do in the sequel is to give a precise description of all
factorisations of a given polynomial (in several variables) with nonnegative
integer coefficients. It turns out that this can be done variable by variable,
and (for a univariate polynomial) that it is useful to classify the polynomials
by number of terms (i.e., evaluation of the polynomial at $1$) rather than
degree. We will prove results for the cases where the number of terms is
at most $10$, or equal to a prime number, with the aid of computer algebra
systems such as MAGMA \cite{Magma}. All software used to produce these proofs
is available from the author.

\subsection{General remarks}

By definition of the polynomial ring (even in infinitely many variables),
every given polynomial involves only finitely many variables of $R$; therefore,
also all its possible factorisations in $R$ together comprise only finitely
many variables, if we identify \emph{isomorphic factors} that differ only in
the labelling of the variables. Thus, we may assume we are in the ring
$R_m=\Z[X_1,\ldots,X_m]$.

Let $N_m\subseteq R_m$ be the subsemiring of polynomials in $X_1,\ldots, X_m$
with nonnegative integer coefficients. Every polynomial $P\in N_m$ has a
\emph{very sparse representation} given by
$$
  P = \sum_{i=1}^n X^{\alpha_i}
$$
where we use the shortcut notation $X^{\gamma} = \prod_{j=1}^m X_j^{\gamma_j}$
for a vector $\gamma=(\gamma_1,\ldots,\gamma_m)$ of nonnegative integer
exponents. The vectors $\alpha_i$ may be repeated in order to accommodate
nontrivial coefficients. We immediately have the following results.

\begin{lemma}
  Let $P\in N_m$ be a sum of $n$ monomials (possibly repeated, all with
  coefficient $1$), and suppose it factors in $N_m$ as $P=ST$, where $S$ has
  $s$ monomials and $T$ has $t$. Then $n=st$. If all monomials of $P$ are
  distinct, then so are those of $S$ and $T$.
\end{lemma}
    
\paragraph*{Proof.}
  Suppose $P=ST$ where $S,T$ are in $N_m$. We have 
  $$
    n=P(1,1,\ldots,1)=S(1,1,\ldots,1) T(1,1,\ldots,1)=st.
  $$
  This proves the first claim.

  Now assume all monomials of $P$ are distinct. Then each term of $P$ arises in
  at most one way as the product of a term of $S$ and a term of $T$, and hence
  $S$ and $T$ cannot have repeated terms.
\proofend

\begin{lemma} \label{LemPrime}
  Let $P\in N_m$ be a sum of $n$ monomials, where $n$ is prime. Then $P$ has
  unique factorisation inside $N_m$.
\end{lemma}

\paragraph*{Proof.}
  Obviously, in any factorisation $P=ST$, either $S$ or $T$ is a monomial.  Now
  the greatest monomial factor is just the greatest common divisor of all terms
  of $P$, and all other monomial factors of $P$ divide it. Thus $P$ is the
  product of a sum of $n$ terms, without common factor, that is irreducible in
  $N_m$, and a monomial. But monomials clearly have unique factorisation in
  $N_m$.
\proofend

~

As soon as the number of terms is not prime, the problem becomes more
difficult. In order to keep track of all possible factorisations, it is
important to consider the order of terms in a polynomial as fixed, and to fix
the order in which we form the terms of the product of two polynomials. This
leads to the following definition.

\begin{defi}
  Let $X_1,\ldots,X_m$ be variables, and let $c_1,\ldots,c_t$ be vectors (of
  length $m$) of nonnegative integers. A \emph{factorisation} of the polynomial
  $\sum_{i=1}^t X^{c_i}$, of the form
  \begin{equation}
  \label{EqFacLem}
    \left(X^{a_1}+\ldots+X^{a_r}\right)\left(X^{b_1}+\ldots+X^{b_s}\right) =
    X^{c_1}+\ldots+X^{c_t},
  \end{equation}
  is given by two sequences $(a_1,\ldots,a_r)$ and $(b_1,\ldots,b_s)$ of
  length-$m$ vectors of nonnegative integers and a bijective mapping $\rho :
  \{1,\ldots,r\} \times \{1,\ldots,s\} \rightarrow \{1,\ldots,t\}$ such that
  \begin{equation} \label{EqBij}
    a_i+b_j=c_{\rho(i,j)} \text{ whenever } 1\le i\le r \text{ and } 
    1\le j\le s.
  \end{equation}
\end{defi}

Of course, for a factorisation to exist, we must have $t=rs$. As remarked
before, we may assume that the polynomials to be factored are not divisible by
variables, and hence we may always assume that $\min_i\{c_{ij}\}=0$ for
$j=1,\ldots,m$. We now show that all possible nonunique factorisations in the
semiring $N_1=\Z_{\ge 0}[X]$ are parametrised by the numbers of terms $r$ and
$s$ and (given $r,s$) the bijections $\rho$.

\begin{lemma} \label{LemLin}
  Let $c_1,\ldots,c_t\in \N^m$ such that $\min_i \{c_{ij}\}=0$ for $1\le j \le
  m$, let $r,s$ be positive integers such that $t=rs$, and let $\rho$ be a
  bijection
  $$
    \rho : \{1,\ldots,r\} \times \{1,\ldots,s\} \rightarrow \{1,\ldots,t\}.
  $$
  Then there exists at most one choice of sequences $(a_1,\ldots,a_r)$ and
  $(b_1,\ldots,b_s)$ of elements of $\N^m$ such that $a_i+b_j=c_{\rho(i,j)}$
  for all $i,j$. If such a choice is possible, then we have
  $\min_i\{a_{ij}\}=\min_i\{b_{ij}\}=0$ and $\{a_{ij}\}\cup \{b_{ij}\} 
  \subseteq \{c_{ij}\}$ for $1\le j\le m$.
\end{lemma}

\paragraph*{Proof.}
  We consider only the first components of all vectors. We get a system
  $$
    a_{i1}+b_{j1} = c_{\rho(i,j),1} \qquad 1\le i\le s, 1\le j\le r
  $$
  of linear equations, to be solved in nonnegative integers. Without loss of
  generality, we assume $c_{11}=0$ and $\rho(1,1)=1$. This means that
  $a_{11}+b_{11}=c_{11}=0$, and hence $a_{11}=b_{11}=0$. The other sums
  $a_{11}+b_{j1}$ and $a_{i1}+b_{11}$ all occur as the left hand side of an
  equation, and hence all $a_{i1}$ and $b_{j1}$ are determined and must occur
  among the $c_{i1}$. The remaining equations then define relations that the
  $c_{i1}$ must satisfy among themselves in order for the system to be soluble.
  It follows that the system has at most one solution for a given $\rho$.

  Finally, it is obvious that the systems for the other components are
  completely analogous and independent from each other, and hence must all be
  satisfied independently. It follows that there is at most one solution.
\proofend

\begin{cor} \label{CorRhodet}
  Let $X_1,\ldots,X_m$ be variables. Given $r$, $s$, $(c_1,\ldots,c_t)\subseteq
  \N^m$, and $\rho$, there exists at most one factorisation of the form
  \eqref{EqFacLem}.
\end{cor}

The last observation in this section shows that the problem of parametrising
all factorisations can be immediately reduced to the univariate case.

\begin{lemma} \label{LemNovars}
  Let $X_1,\ldots,X_m$ be variables, and let $c_1,\ldots,c_t\in\N^m$.
  A quintuple $(r,s,\rho,(a_1,\ldots,a_r),(b_1,\ldots,b_s))$ is a
  factorisation of $\sum_{i=1}^t X^{c_i}$ in $N_m$ if and only if
  $(r,s,\rho,(a_{1j},\ldots,a_{rj}),(b_{1j},\ldots,b_{sj}))$ is a 
  factorisation of $\sum_{i=1}^t X_j^{c_{ij}}$ inside $\N[X_j]$, for
  $j=1,\ldots,m$.
\end{lemma}

\paragraph*{Proof.}
  A polynomial in $N_m$ in very sparse representation is determined by its
  exponent vectors $c_1,\ldots,c_t$. If we set all variables except $X_1$, say,
  to one, this corresponds to setting all components of the $c_i$, except the
  first, to zero. If we do this for all the $X_j$ separately, we see that no
  information in the $c_i$ is lost; in fact, we map the matrix $(c_{ij})$ to
  the set of its columns. In particular, factorisations are preserved, and
  factorisations with respect to each variable separately can be pieced
  together to obtain a factorisation of the original polynomial.
\proofend

Note that the maps described in the proof are not the same as the semiring
homomorphisms of evaluating all variables (except one) at $1$. For example, the
polynomial $XY+1$ in very sparse representation maps to $(X+1,Y+1)$ as usual,
but $X+Y$ maps to $(X+1,1+Y)$, which is different since the order of terms is
fixed. This allows us to reconstruct polynomials in $N_m$ from their images
under these partial evaluations; in fact, the number of terms (in the very
sparse sense) does not change under our maps.

\subsection{The univariate case}

From now on, we restrict attention to the case of polynomials in $N_1$. More
precisely, let $t$ be a nonnegative integer; we consider polynomials $P\in N_1$
such that $P(1)=t$ and $P(0)>0$, so that $P$ is not divisible by the variable.
By what we just proved, it follows that we can obtain all possible
factorisations of $P$ inside $N_1$ by listing all nontrivial factorisations
$t=rs$ (i.e., $r>1$ and $s>1$), and for every pair $(r,s)$, listing all
possible bijections $\rho:\{1,\ldots,r\}\times \{1,\ldots,s\} \rightarrow
\{1,\ldots,t\}$. To complete the arguments, we will need to prove that the
factors thus found are irreducible; if they are not, we will have to consider
products of more than two factors. We will give complete descriptions of the
cases $t=1,\ldots,12$; of course, the cases where $t$ is prime were already
dealt with by Lemma \ref{LemPrime}.

Suppose $r$ and $s$ are given. We will assume that factorisations have the form
\eqref{EqFacLem}, where moreover we have
\begin{gather} 
  \label{EqSorted1} 0=c_1 \le c_2 \le \ldots \le c_t; \\
  \label{EqSorted2}
  0=a_1 \le \ldots \le a_r \qquad \text{and} \qquad 0=b_1 \le \ldots \le b_s.
\end{gather}
As in the proof of Lemma \ref{LemLin}, the exponents $a_1+b_j=b_j$ correspond
to $c_{\rho(1,j)}$, and $a_i+b_1=a_i$ to $c_{\rho(i,1)}$. Therefore, we obtain
the equations
\begin{equation} \label{EqImplicit}
  c_{\rho(i,j)} = c_{\rho(i,1)} + c_{\rho(1,j)} 
     \qquad (2\le i\le r,2\le j\le s).
\end{equation}
Thus, for every bijection $\rho$, the polynomials possessing a factorisation of
the form \eqref{EqFacLem} with bijection $\rho$ correspond to the solutions in
the (nonnegative integer) variables $c_i$ of the system of linear equations and
inequalities made up by \eqref{EqSorted1} together with \eqref{EqImplicit}.
Such systems are closely related to \emph{integer linear programs}. The
difference is that there is no cost function to be optimised; instead, we are
interested in \emph{all solutions}, which hopefully will admit a compact
description. If the polynomial $P$ is to have two different factorisations,
then the $c_i$ must satisfy two such linear programs simultaneously, namely
those corresponding to two different $\rho$.

To facilitate the analysis, we can use the fact that \eqref{EqImplicit} also
implies a partial ordering of the $c_i$; in fact, because all exponents are
nonnegative, we have
\begin{equation} \label{EqPartialWeak}
  c_{\rho(i,j)} \ge c_{\rho(i,1)} \text{ and } c_{\rho(i,j)} \ge c_{\rho(1,j)} 
     \qquad (2\le i\le r,\, 2\le j\le s).
\end{equation}
Furthermore, in the end we will drop the fixed term ordering of the factors,
so that we may assume \eqref{EqSorted2}. The implication is that
\begin{equation} \label{EqPartialStrong}
  c_{\rho(i,j)} \ge c_{\rho(i',j')} 
     \qquad (1\le i' \le i\le r,\, 1\le j'\le j\le s).
\end{equation}
The number of pairs of $c_i$ whose ordering is thus determined is
$$
  \binom{r}{2}\binom{s}{2} + r\binom{s}{2} + \binom{r}{2}s = 
  \frac{r^2s^2 + rs^2 + r^2s -3rs}4,
$$
to be compared with the total number $\binom{rs}{2}=\frac{r^2s^2-rs}2$ of
pairs. The idea is that the partial ordering \eqref{EqPartialStrong} severely
limits the possibilities for $\rho$.

In the (computer-aided) proofs below, we will repeatedly compute all
possibilities for the bijection $\rho$ that are compatible with the partial
order \eqref{EqPartialStrong}. It is readily seen that this problem is
equivalent to the \emph{topological sorting} of a given directed acyclic graph
(dag) with nodes $\{1,\ldots,rs\}$.  Here the graph represents the few ordering
relationships that are known, and the task is to produce a total ordering of
the nodes such that node $i$ comes before node $j$ whenever the graph contains
an edge $(i,j)$. Such orderings exist if and only if the graph is acyclic, and
are not unique. To find all possible such orderings, we apply Kahn's algorithm,
with backtracking. (It is possible to use a depth-first search to find a
topological ordering, and this is the algorithm given in \cite[Section
23.4]{CormenLeisersonRivest:1990}, but this algorithm is unable to derive
\emph{all} topological orderings --- cf.\ Exercise 23.4-2 in
\cite{CormenLeisersonRivest:1990}. In fact, on the graph with three nodes
and edges $(1,3)$ and $(2,3)$, depth-first search will never produce the
ordering $123$, which is perfectly valid.)

Kahn's algorithm works as follows: take any node with zero in-degree, write it
down, remove it and all its incident edges from the graph, and continue
recursively with the remainder until the graph is trivial. By backtracking on
every choice that we make, we smoothly produce all possible topological
orderings.

\subsection{Results for polynomials with few terms}

Using the above framework, we arrive at the following results.

\begin{theorem} \label{Thm4}
  Let $P\in N_1$ be a sum of $4$ monomials, not necessarily distinct. Then $P$
  has unique factorisation inside $N_1$.
\end{theorem}

\paragraph*{Proof.}
  We may assume the variable does not divide $P$. Putting $t=4$, the only
  nontrivial factorisation is one having $r=s=2$. 
  
  Let us see which bijections $\rho$ can occur in the factorisation
  $$
    ( 1 + X^{a_2} ) ( 1 + X^{b_2} ) = 1 + X^{c_2} + X^{c_3} + X^{c_4}.
  $$
  We may put $\rho(1,1)=1$ and $\rho(2,2)=4$, because $c_{\rho(2,2)}=a_2+b_2$
  must be maximal among the $c_i$. The only remaining choice for $\rho$ is
  whether $\rho(1,2)=2$ or $3$ (i.e., whether $c_2=a_2$ or $c_3=a_2$). But this
  choice corresponds to exchanging the factors on the left, and this leaves the
  factorisation is essentially unchanged. Thus $\rho$ is fixed, and $P$ is
  either irreducible (if $c_4\ne c_2+c+3$) or breaks up uniquely into two
  factors of $2$ terms each, by Corollary \ref{CorRhodet}.
\proofend

\begin{theorem} \label{Thm6}
  Let $P\in N_1$ be a sum of $6$ monomials, not necessarily distinct. If $P$
  has non-unique factorisation, then it is of the form
  $$
    X^a(1 + X^b + X^{2b} + X^{3b} + X^{4b} + X^{5b})
  $$
  for nonnegative integers $a$ and $b\ge 1$.
\end{theorem}

Note that $1+X+X^2+X^3+X^4+X^5= (1+X^2+X^4)(1+X)= (1+X+X^2)(1+X^3)$. All 
factors are irreducible in $N_1$ by Lemma \ref{LemPrime}.

\paragraph*{Proof.}
  We may assume the variable does not divide $P$. Any nontrivial factorisation
  is of the form
  $$
    (1 + X^{a_2} + X^{a_3})(1+X^{b_2}).
  $$
  We use the topological sorting algorithm described above to find all
  bijections $\rho$ that are compatible with the partial ordering
  \eqref{EqPartialStrong}, and find that there are the following $5$
  possibilities (to be read such that the first pair corresponds to $c_1$, the
  next to $c_2$, and so on):
  $$
    \begin{array}{c|c|c|c|c|c|c}
    \text{Case:} & c_1 & c_2 & c_3 & c_4 & c_5 & c_6 \\ \hline
    1 & (1,1) & (1,2) & (2,1) & (2,2) & (3,1) & (3,2) \\
    2 & (1,1) & (1,2) & (2,1) & (3,1) & (2,2) & (3,2) \\
    3 & (1,1) & (2,1) & (1,2) & (2,2) & (3,1) & (3,2) \\
    4 & (1,1) & (2,1) & (1,2) & (3,1) & (2,2) & (3,2) \\
    5 & (1,1) & (2,1) & (3,1) & (1,2) & (2,2) & (3,2) \\ \hline
    \end{array} 
  $$
  We now check whether any combination of two $\rho$ admits a simultaneous
  solution of the respective sets of equations \eqref{EqImplicit}. For example,
  the first line implies the equations
  $$
    c_4 = c_2+c_3 \text{ and } c_6=c_2+c_5,
  $$
  whereas the fifth line gives
  $$
    c_5=c_2+c_4 \text{ and } c_6=c_3+c_4.
  $$
  The general solution to these four equations combined is
  $$
    c_2 \text{ free}, \; c_3=2c_2, \; c_4=3c_2, \; c_5=4c_2, \; c_6=5c_2,
  $$
  and this corresponds to the factorisation given in the Theorem. To see this,
  note that the bijection $\rho$ from the first line results in
  $a_2=c_3$, $a_3=c_5$, and $b_2=c_2$, whereas the fifth line gives
  $a_2=c_2$, $a_3=c_3$, and $b_2=c_4$.

  It is easily (but tediously) verified that in all the factorisations given
  by combining two other lines of the table, the left and right hand sides
  coincide. Hence the only nonunique factorisation is the one given in the
  Theorem, and {\it a fortiori}, polynomials that possess three inequivalent
  factorisations do not exist.
\proofend

\begin{theorem} \label{Thm8}
  Let $P\in N_1$ be a sum of $8$ monomials, not necessarily distinct. Then $P$
  has unique factorisation inside $N_1$.
\end{theorem}

\paragraph*{Proof.}
  We may assume the variable does not divide $P$. If $P$ is reducible in $N_1$,
  then it can be written as a product of factors with $2$ and $4$ terms, where
  the latter possibly again breaks up as a product of factors with $2$ terms
  each. We will first establish all possible nonunique factorisations of type
  $$
    (1 + X^{a_2})(1 + X^{b_2} + X^{b_3} + X^{b_4}).
  $$
  All bijections $\rho:\{1,2\}\times\{1,2,3,4\}\rightarrow\{1,\ldots,8\}$
  that are compatible with the partial ordering \eqref{EqPartialStrong} are
  enumerated by MAGMA as
  $$
    \begin{array}{c|c|c|c|c|c|c|c|c}
    \text{Case:} & c_1 & c_2 & c_3 & c_4 & c_5 & c_6 & c_7 & c_8 \\ \hline
    1 & (1,1) & (1,2) & (1,3) & (1,4) & (2,1) & (2,2) & (2,3) & (2,4) \\
    2 & (1,1) & (1,2) & (1,3) & (2,1) & (1,4) & (2,2) & (2,3) & (2,4) \\
    3 & (1,1) & (1,2) & (1,3) & (2,1) & (2,2) & (1,4) & (2,3) & (2,4) \\
    4 & (1,1) & (1,2) & (1,3) & (2,1) & (2,2) & (2,3) & (1,4) & (2,4) \\
    5 & (1,1) & (1,2) & (2,1) & (1,3) & (1,4) & (2,2) & (2,3) & (2,4) \\
    6 & (1,1) & (1,2) & (2,1) & (1,3) & (2,2) & (1,4) & (2,3) & (2,4) \\
    7 & (1,1) & (1,2) & (2,1) & (1,3) & (2,2) & (2,3) & (1,4) & (2,4) \\
    8 & (1,1) & (1,2) & (2,1) & (2,2) & (1,3) & (1,4) & (2,3) & (2,4) \\
    9 & (1,1) & (1,2) & (2,1) & (2,2) & (1,3) & (2,3) & (1,4) & (2,4) \\
    10& (1,1) & (2,1) & (1,2) & (1,3) & (1,4) & (2,2) & (2,3) & (2,4) \\
    11& (1,1) & (2,1) & (1,2) & (1,3) & (2,2) & (1,4) & (2,3) & (2,4) \\
    12& (1,1) & (2,1) & (1,2) & (1,3) & (2,2) & (2,3) & (1,4) & (2,4) \\
    13& (1,1) & (2,1) & (1,2) & (2,2) & (1,3) & (1,4) & (2,3) & (2,4) \\
    14& (1,1) & (2,1) & (1,2) & (2,2) & (1,3) & (2,3) & (1,4) & (2,4) \\ \hline
    \end{array} 
  $$
  The equations \eqref{EqImplicit} are
  \begin{align*}
    c_{\rho(2,2)} &= c_{\rho(2,1)} + c_{\rho(1,2)} = a_2 + b_2; \\
    c_{\rho(2,3)} &= c_{\rho(2,1)} + c_{\rho(1,3)} = a_2 + b_3; \\
    c_{\rho(2,4)} &= c_{\rho(2,1)} + c_{\rho(1,4)} = a_2 + b_4.
  \end{align*}
  As a typical example, we will consider polynomials admitting factorisations
  under both bijections $1$ and $8$. Besides $c_1=0$, which always holds, we
  find
  \begin{gather*}
    c_6 = c_2 + c_5; \qquad c_7 = c_3 + c_5; \qquad c_8 = c_4 + c_5; \\
    c_4 = c_2 + c_3; \qquad c_7 = c_3 + c_5; \qquad c_8 = c_3 + c_6,
  \end{gather*}
  which is obviously solvable taking $c_2$, $c_3$, and $c_5$ as free
  (nonnegative integer) parameters. The corresponding inequivalent
  factorisations are
  $$
    ( 1 + X^{c_5} ) ( 1 + X^{c_2} + X^{c_3} + X^{c_2+c_3} ) = 
    ( 1 + X^{c_3} ) ( 1 + X^{c_2} + X^{c_5} + X^{c_2+c_5} ).
  $$
  However, in fact both can be further reduced to
  $$
    ( 1 + X^{c_2} ) (1 + X^{c_3} ) ( 1 + X^{c_5} ).
  $$
  Automated verification by MAGMA shows that likewise all pairs of inequivalent
  $2\times 4$-factorisations that arise by combining two cases from the table
  can be further reduced to equivalent $2\times 2\times 2$-factorisations. This
  proves the theorem.
\proofend

\begin{theorem} \label{Thm9}
  Let $P\in N_1$ be a sum of $9$ monomials, not necessarily distinct. Then $P$
  has unique factorisation inside $N_1$.
\end{theorem}

\paragraph*{Proof.}
  We may assume the variable does not divide $P$. Putting $t=9$, the only
  nontrivial factorisation is one having $r=s=3$. 
  
  We establish all factorisations of the form
  $$
    ( 1 + X^{a_2} + X^{a_3} ) ( 1 + X^{b_2} + X^{b_3} ) = 
      1 + X^{c_2} + X^{c_3} + X^{c_4} + X^{c_5} + X^{c_6} + X^{c_7} +
      X^{c_8} + X^{c_9}.
  $$
  Exchanging the factors does not change the factorisation, so in addition to
  the partial ordering \eqref{EqPartialStrong}, we also require that $a_2 \ge
  b_2$, hence $c_{\rho(2,1)} \ge c_{\rho(1,2)}$. All bijections $\rho$ that
  are compatible with these orderings are enumerated by MAGMA to be
  $$
    \begin{array}{c|c|c|c|c|c|c|c|c|c}
    \text{Case:} & c_1 & c_2 & c_3 & c_4 & c_5 & c_6 & c_7 & c_8 & c_9 \\ \hline
      1 & (1,1) & (1,2) & (1,3) & (2,1) & (2,2) & (2,3) & (3,1) & (3,2) & (3,3) \\
      2 & (1,1) & (1,2) & (1,3) & (2,1) & (2,2) & (3,1) & (2,3) & (3,2) & (3,3) \\
      3 & (1,1) & (1,2) & (1,3) & (2,1) & (2,2) & (3,1) & (3,2) & (2,3) & (3,3) \\
      4 & (1,1) & (1,2) & (1,3) & (2,1) & (3,1) & (2,2) & (2,3) & (3,2) & (3,3) \\
      5 & (1,1) & (1,2) & (1,3) & (2,1) & (3,1) & (2,2) & (3,2) & (2,3) & (3,3) \\
      6 & (1,1) & (1,2) & (2,1) & (1,3) & (2,2) & (2,3) & (3,1) & (3,2) & (3,3) \\
      7 & (1,1) & (1,2) & (2,1) & (1,3) & (2,2) & (3,1) & (2,3) & (3,2) & (3,3) \\
      8 & (1,1) & (1,2) & (2,1) & (1,3) & (2,2) & (3,1) & (3,2) & (2,3) & (3,3) \\
      9 & (1,1) & (1,2) & (2,1) & (1,3) & (3,1) & (2,2) & (2,3) & (3,2) & (3,3) \\
     10 & (1,1) & (1,2) & (2,1) & (1,3) & (3,1) & (2,2) & (3,2) & (2,3) & (3,3) \\
     11 & (1,1) & (1,2) & (2,1) & (2,2) & (1,3) & (2,3) & (3,1) & (3,2) & (3,3) \\
     12 & (1,1) & (1,2) & (2,1) & (2,2) & (1,3) & (3,1) & (2,3) & (3,2) & (3,3) \\
     13 & (1,1) & (1,2) & (2,1) & (2,2) & (1,3) & (3,1) & (3,2) & (2,3) & (3,3) \\
     14 & (1,1) & (1,2) & (2,1) & (2,2) & (3,1) & (1,3) & (2,3) & (3,2) & (3,3) \\
     15 & (1,1) & (1,2) & (2,1) & (2,2) & (3,1) & (1,3) & (3,2) & (2,3) & (3,3) \\
     16 & (1,1) & (1,2) & (2,1) & (2,2) & (3,1) & (3,2) & (1,3) & (2,3) & (3,3) \\
     17 & (1,1) & (1,2) & (2,1) & (3,1) & (1,3) & (2,2) & (2,3) & (3,2) & (3,3) \\
     18 & (1,1) & (1,2) & (2,1) & (3,1) & (1,3) & (2,2) & (3,2) & (2,3) & (3,3) \\
     19 & (1,1) & (1,2) & (2,1) & (3,1) & (2,2) & (1,3) & (2,3) & (3,2) & (3,3) \\
     20 & (1,1) & (1,2) & (2,1) & (3,1) & (2,2) & (1,3) & (3,2) & (2,3) & (3,3) \\
     21 & (1,1) & (1,2) & (2,1) & (3,1) & (2,2) & (3,2) & (1,3) & (2,3) & (3,3) 
\\
     \hline
    \end{array} 
  $$
  As before, we consider the systems of $8$ equations in $8$ variables (putting
  the trivial $c_1=0$ apart) that arise by combining two possibilities for
  $\rho$ as given in the table. Among the $8$ variables, at most four are
  conceptually suited as free parameters, namely those corresponding to one of
  the variables $a_2$, $a_3$, $b_2$, and $b_3$ \emph{under both bijections
  simultaneously}. It is interesting to note that it is always possible to
  choose the free parameters among this small set of ``interesting'' variables
  (it would be even more interesting to prove this beforehand), even if some of
  the occurring systems have rank as small as $5$ (out of $8$).

  We let MAGMA run through all combinations of two bijections, and the result
  is that in all cases, the obtained factorisations are equivalent.
\proofend

\begin{theorem} \label{Thm10}
  Let $P\in N_1$ be a sum of $10$ monomials, not necessarily distinct. If $P$
  has non-unique factorisation, then it is of the form
  \begin{align*}
    (1+X)(1+X^2+X^4+X^6+X^8) &= (1+X^5)(1+X+X^2+X^3+X^4); \\
    (1+X)(1+X^4+X^6+X^8+X^{12}) &= (1+X^5)(1+X+X^4+X^7+X^8); \\ 
    (1+X^3)(1+X^2+X^4+X^6+X^8) &= (1+X^5)(1+X^2+X^3+X^4+X^6),
  \end{align*}
  up to multiplication by a power of $X$, and up to replacing $X$ by a power,
  or it has one of the forms
  \begin{align*}
    (1 + X^{3b}) ( 1 + X^a + X^b + X^{a+b} + X^{a+2b} ) &=
      ( 1 + X^b ) ( 1 + X^a + X^{a+2b} + X^{3b} + X^{a+4b} ); \\
    (1 + X^{3b}) ( 1 + X^a + X^b + X^{a+b} + X^{2b} ) &=
      ( 1 + X^b ) ( 1 + X^a + X^{2b} + X^{a+3b} + X^{4b} )
  \end{align*}
  for integers $a\ge 0$ and $b\ge 1$, up to multiplication by a power of $X$.
\end{theorem}

The products in the Theorem evaluate to
\begin{align*}
    & 1 + X + X^2 + X^3 + X^4 + X^5 + X^6 + X^7 + X^8 + X^9; \\
    & 1 + X + X^4 + X^5 + X^6 + X^7 + X^8 + X^9 + X^{12} + X^{13}; \\
    & 1 + X^2 + X^3 + X^4 + X^5 + X^6 + X^7 + X^8 + X^9 + X^{11}; \\
    & 1 + X^a + X^b + X^{a+b} + X^{a+2b} + X^{3b} + X^{a+3b} + X^{4b}
      + X^{a+4b} + X^{5b}; \\
    & 1 + X^a + X^b + X^{a+b} + X^{2b} + X^{3b} + X^{a+3b} + X^{4b}
      + X^{a+4b} + X^{5b},
\end{align*}
respectively. All factors are irreducible in $N_1$ by Lemma \ref{LemPrime}.
Note that the first three cases are all self-reciprocal; the reciprocals of all
polynomials in the two $2$-parameter families again belong to one of these two
families.

\paragraph*{Proof.}
  We may assume the variable does not divide $P$. Any nontrivial factorisation
  is of the form
  $$
    (1 + X^{a_2}) (1 + X^{b_2} + X^{b_3} + X^{b_4} + X^{b_5} ).
  $$
  We use topological sorting and backtracking to find all bijections $\rho$
  that are compatible with the partial ordering \eqref{EqPartialStrong}, and
  find the following cases:
  {\small
  $$
    \begin{array}{c|c|c|c|c|c|c|c|c|c|c}
    \text{Case:} & c_1 & c_2 & c_3 & c_4 & c_5 & c_6 & c_7 & c_8 & c_9 & c_{10} \\ \hline
      1 & (1,1) & (1,2) & (1,3) & (1,4) & (1,5) & (2,1) & (2,2) & (2,3) & (2,4) & (2,5) \\
      2 & (1,1) & (1,2) & (1,3) & (1,4) & (2,1) & (1,5) & (2,2) & (2,3) & (2,4) & (2,5) \\
      3 & (1,1) & (1,2) & (1,3) & (1,4) & (2,1) & (2,2) & (1,5) & (2,3) & (2,4) & (2,5) \\
      4 & (1,1) & (1,2) & (1,3) & (1,4) & (2,1) & (2,2) & (2,3) & (1,5) & (2,4) & (2,5) \\
      5 & (1,1) & (1,2) & (1,3) & (1,4) & (2,1) & (2,2) & (2,3) & (2,4) & (1,5) & (2,5) \\
      6 & (1,1) & (1,2) & (1,3) & (2,1) & (1,4) & (1,5) & (2,2) & (2,3) & (2,4) & (2,5) \\
      7 & (1,1) & (1,2) & (1,3) & (2,1) & (1,4) & (2,2) & (1,5) & (2,3) & (2,4) & (2,5) \\
      8 & (1,1) & (1,2) & (1,3) & (2,1) & (1,4) & (2,2) & (2,3) & (1,5) & (2,4) & (2,5) \\
      9 & (1,1) & (1,2) & (1,3) & (2,1) & (1,4) & (2,2) & (2,3) & (2,4) & (1,5) & (2,5) \\
     10 & (1,1) & (1,2) & (1,3) & (2,1) & (2,2) & (1,4) & (1,5) & (2,3) & (2,4) & (2,5) \\
     11 & (1,1) & (1,2) & (1,3) & (2,1) & (2,2) & (1,4) & (2,3) & (1,5) & (2,4) & (2,5) \\
     12 & (1,1) & (1,2) & (1,3) & (2,1) & (2,2) & (1,4) & (2,3) & (2,4) & (1,5) & (2,5) \\
     13 & (1,1) & (1,2) & (1,3) & (2,1) & (2,2) & (2,3) & (1,4) & (1,5) & (2,4) & (2,5) \\
     14 & (1,1) & (1,2) & (1,3) & (2,1) & (2,2) & (2,3) & (1,4) & (2,4) & (1,5) & (2,5) \\
     15 & (1,1) & (1,2) & (2,1) & (1,3) & (1,4) & (1,5) & (2,2) & (2,3) & (2,4) & (2,5) \\
     16 & (1,1) & (1,2) & (2,1) & (1,3) & (1,4) & (2,2) & (1,5) & (2,3) & (2,4) & (2,5) \\
     17 & (1,1) & (1,2) & (2,1) & (1,3) & (1,4) & (2,2) & (2,3) & (1,5) & (2,4) & (2,5) \\
     18 & (1,1) & (1,2) & (2,1) & (1,3) & (1,4) & (2,2) & (2,3) & (2,4) & (1,5) & (2,5) \\
     19 & (1,1) & (1,2) & (2,1) & (1,3) & (2,2) & (1,4) & (1,5) & (2,3) & (2,4) & (2,5) \\
     20 & (1,1) & (1,2) & (2,1) & (1,3) & (2,2) & (1,4) & (2,3) & (1,5) & (2,4) & (2,5) \\
     21 & (1,1) & (1,2) & (2,1) & (1,3) & (2,2) & (1,4) & (2,3) & (2,4) & (1,5) & (2,5) \\
     22 & (1,1) & (1,2) & (2,1) & (1,3) & (2,2) & (2,3) & (1,4) & (1,5) & (2,4) & (2,5) \\
     23 & (1,1) & (1,2) & (2,1) & (1,3) & (2,2) & (2,3) & (1,4) & (2,4) & (1,5) & (2,5) \\
     24 & (1,1) & (1,2) & (2,1) & (2,2) & (1,3) & (1,4) & (1,5) & (2,3) & (2,4) & (2,5) \\
     25 & (1,1) & (1,2) & (2,1) & (2,2) & (1,3) & (1,4) & (2,3) & (1,5) & (2,4) & (2,5) \\
     26 & (1,1) & (1,2) & (2,1) & (2,2) & (1,3) & (1,4) & (2,3) & (2,4) & (1,5) & (2,5) \\
     27 & (1,1) & (1,2) & (2,1) & (2,2) & (1,3) & (2,3) & (1,4) & (1,5) & (2,4) & (2,5) \\
     28 & (1,1) & (1,2) & (2,1) & (2,2) & (1,3) & (2,3) & (1,4) & (2,4) & (1,5) & (2,5) \\
     29 & (1,1) & (2,1) & (1,2) & (1,3) & (1,4) & (1,5) & (2,2) & (2,3) & (2,4) & (2,5) \\
     30 & (1,1) & (2,1) & (1,2) & (1,3) & (1,4) & (2,2) & (1,5) & (2,3) & (2,4) & (2,5) \\
     31 & (1,1) & (2,1) & (1,2) & (1,3) & (1,4) & (2,2) & (2,3) & (1,5) & (2,4) & (2,5) \\
     32 & (1,1) & (2,1) & (1,2) & (1,3) & (1,4) & (2,2) & (2,3) & (2,4) & (1,5) & (2,5) \\
     33 & (1,1) & (2,1) & (1,2) & (1,3) & (2,2) & (1,4) & (1,5) & (2,3) & (2,4) & (2,5) \\
     34 & (1,1) & (2,1) & (1,2) & (1,3) & (2,2) & (1,4) & (2,3) & (1,5) & (2,4) & (2,5) \\
     35 & (1,1) & (2,1) & (1,2) & (1,3) & (2,2) & (1,4) & (2,3) & (2,4) & (1,5) & (2,5) \\
     36 & (1,1) & (2,1) & (1,2) & (1,3) & (2,2) & (2,3) & (1,4) & (1,5) & (2,4) & (2,5) \\
     37 & (1,1) & (2,1) & (1,2) & (1,3) & (2,2) & (2,3) & (1,4) & (2,4) & (1,5) & (2,5) \\
     38 & (1,1) & (2,1) & (1,2) & (2,2) & (1,3) & (1,4) & (1,5) & (2,3) & (2,4) & (2,5) \\
     39 & (1,1) & (2,1) & (1,2) & (2,2) & (1,3) & (1,4) & (2,3) & (1,5) & (2,4) & (2,5) \\
     40 & (1,1) & (2,1) & (1,2) & (2,2) & (1,3) & (1,4) & (2,3) & (2,4) & (1,5) & (2,5) \\
     41 & (1,1) & (2,1) & (1,2) & (2,2) & (1,3) & (2,3) & (1,4) & (1,5) & (2,4) & (2,5) \\
     42 & (1,1) & (2,1) & (1,2) & (2,2) & (1,3) & (2,3) & (1,4) & (2,4) & (1,5) & (2,5) \\
     \hline
    \end{array} 
  $$
  }
  Solving the linear systems obtained by every combination of two entries in 
  the table, this time we find $102$ pairs with nonequivalent factorisations.
  A good way to classify the cases is to look at the two values for the
  variable $a_2$, i.e., the nontrivial exponent of $X$ in the left-hand factor.

  In some cases, such as the combination Cases $5$ and $38$, the two choices
  for the left-hand factor are $(1+X^a)$ and $(1+X^{-a})$, for a free parameter
  $a$. As the parameter must be integer and nonnegative, it is forced to be
  $0$, and the corresponding factorisations are equivalent after all.

  In the combinations $1$ and $42$, as well as $10$ and $42$, the left-hand
  factors are $(1+X^a)$ and $(1+X^{5a})$, giving the first and second case in
  the Theorem.

  In the combinations $(1,9)$, $(1,17)$, $(1,30)$, $(2,12)$, $(2,20)$, 
  $(2,33)$, $(3,14)$, $(3,22)$, $(6,25)$, $(6,38)$, and $(7,27)$, the 
  left-hand factors are $(1+X^{3a})$ and $(1+X^{5a})$, giving the third case in
  the Theorem (all combinations yield the same right-hand factors).

  In all other combinations, the ratio of the exponents of $X$ in the left-hand
  factors is $1:3$, and we find that all cases are covered by the 
  two-parameter families given in the Theorem. Here, some combinations give
  one of the two-parameter families, while others give a specialisation of
  either family, obtained by putting $a=kb$ for $k\in\{0,1,2,3,4,5\}$.
\proofend

In fact (an observation of Wilfried Imrich), each of the theorems given above
yields a slightly more general result. Namely, if instead of considering just
terms with coefficient $1$, we give each term of our polynomials a positive
integer coefficient, then the same proofs show that there is at most one 
factorisation of $P$ of the given form (for example, a factorisation of a
$10$-term polynomial into a product of a $2$-term and a $5$-term polynomial). 
For it to hold, now also the coefficients must satisfy some easy relations.
However, if we have nontrivial coefficients, the number of terms a factor
has need not divide the number of terms of the polynomial itself, and therefore
(unfortunately) this is not enough to show that all quadrinomials have unique
factorisation in $N_1$.

The method above clearly reduces the classification of all non-unique 
factorisations of univariate polynomials with a given number of terms $n$ to a
finite computation. However, it is also clear that the cost of these
computations increases exponentially with $n$. Hence we refrain from attacking
the complicated case $n=12$. We only note that factorisation of $12$-term
polynomials in $N_1$ is non-unique, as witnessed by N\"usken's example
\eqref{EqNusken}. Also, although cases where $n$ is a prime power apparently
tend to possess few nonunique factorisations or none at all, there do exist
examples, such as
\begin{equation} \label{Eq16}
  (1+X^3+X^5+X^6)(1+X+X^2+X^4) = (1+X)(1+X^2)(1+2X^4+X^7)
\end{equation}
for $n=16$, which also shows that the monoid $N_1$ is not half-factorial.

\section*{Acknowledgement}

I would like to thank Richard Hammack and Wilfried Imrich for introducing me
to the subject in such an agreeable way.

\def\polhk#1{\setbox0=\hbox{#1}{\ooalign{\hidewidth
  \lower1.5ex\hbox{`}\hidewidth\crcr\unhbox0}}} \def\cprime{$'$}
  \hyphenation{ma-the-ma-ti-ques}

\end{document}